\documentclass[11pt]{article}
\usepackage[english]{babel}
\usepackage{times}
\usepackage[T1]{fontenc}
\usepackage{amsmath, amsthm, amsfonts}
\usepackage{amssymb}
\usepackage{amscd}
\usepackage{verbatim}

\title{ Trace Hardy--Sobolev--Mazy'a  inequalities \\  for the half
fractional Laplacian}

\author{\Large Stathis Filippas$^{1,4}$,~~  Luisa Moschini$^{2}$~
\&~ Achilles Tertikas$^{3,4}$  \\
                                                                           \\
        Department of Applied Mathematics$^{1}$ \\
         University of Crete,
         71409 Heraklion,  Greece \\
        filippas@tem.uoc.gr\\
                                          \\
                                          Dipartimento di Scienze di Base  ed Applicate per l'Ingegneria $^2$,\\
 University of Rome ``La Sapienza'', 00185 Rome, Italy \\
luisa.moschini@sbai.uniroma1.it \\
\\
 Department of Mathematics$^{3}$ \\
         University of Crete,
         71409 Heraklion,  Greece \\
          tertikas@math.uoc.gr\\
                                     \\
        Institute of Applied and Computational Mathematics$^4$, \\
        FORTH, 71110 Heraklion, Greece \\
    \\ }

\begin{document}
\date{}
\newcommand{\ana}{\nabla}
\newcommand{\R}{I \!  \! R}
\newcommand{\N}{I \!  \! N}
\newcommand{\Ren}{ I \! \! R^n}
\maketitle


\newcommand{\be}{\begin{equation}}
\newcommand{\ee}{\end{equation}}
\newcommand{\bea}{\begin{eqnarray}}
\newcommand{\eea}{\end{eqnarray}}
\newcommand{\bean}{\begin{eqnarray*}}
\newcommand{\eean}{\end{eqnarray*}}
\newcommand{\la}{\label}
\newcommand{\bx}{\bar{x}}


\newcommand{\xa}{\alpha}
\newcommand{\xb}{\beta}
\newcommand{\xg}{\gamma}
\newcommand{\xG}{\Gamma}
\newcommand{\xd}{\delta}
\newcommand{\xD}{\Delta}
\newcommand{\xe}{\varepsilon}
\newcommand{\xz}{\zeta}
\newcommand{\xh}{\eta}
\newcommand{\Th}{\Theta}
\newcommand{\xk}{\kappa}
\newcommand{\xl}{\lambda}
\newcommand{\xL}{\Lambda}
\newcommand{\CC}{\mathcal{C}}
\newcommand{\xr}{\rho}
\newcommand{\xs}{\sigma}
\newcommand{\xS}{\Sigma}
\newcommand{\xo}{\omega}
\newcommand{\xO}{\Omega}
\newcommand{\tA}{\tilde{A}}
\newcommand{\tC}{\bar{k}_s}
\newcommand{\finedim}{{\hfill $\Box$}}

\newcommand{\dino}{\int_0^{+\infty}\int_{\xO}}
\newcommand{\dinoC}{\int_0^{+\infty}\int_{\CC \xO}}
\newcommand{\dinoc}{\int_0^{+\infty}\int_{\xO^c}}
\newcommand{\dinr}{\int_0^{+\infty}\int_{\Ren}}
\newcommand{\dd}{dxdy}
\newcommand{\ra}{\rightarrow}
\newcommand{\rft}{\rightarrow +\infty}
\newcounter{newsection}


\newtheorem{theorem}{Theorem}[section]
\newtheorem{lemma}[theorem]{Lemma}
\newtheorem{prop}[theorem]{Proposition}
\newtheorem{coro}[theorem]{Corollary}
\newtheorem{defin}[theorem]{Definition}
\newcounter{newsec} \renewcommand{\theequation}{\thesection.\arabic{equation}}

\begin{abstract}
In this work we  establish  trace
Hardy-Sobolev-Maz'ya inequalities with best Hardy constants,
 for weakly  mean convex  domains. We accomplish this by obtaining
 a new weighted  Hardy type estimate which is of independent inerest. We then  produce
Hardy-Sobolev-Maz'ya inequalities for the spectral half Laplacian. This covers a critical
case left open in \cite{FMT1}.
\end{abstract}




\setcounter{equation}{0}
\section{Introduction}\la{intro}

The  Hardy-Sobolev-Maz'ya  (HSM)  inequalities combine the   Sobolev and   Hardy  terms,   the latter with best constant.
 For instance, for  the regular (local) Laplacian and for a domain  $\xO \subsetneqq \Ren$
   $n \geq 3$ it states that, if
$d(x)= {\rm dist}(x, \partial \xO)$,
 there exists a positive
constant $c$ such that
\begin{equation} \la{01}
\int_{\xO}|\nabla u|^2dx \geq \frac{1}{4}
\int_{\xO}\frac{|u|^2}{d^2(x)}dx + c \left(
\int_{\xO}|u|^{\frac{2n}{n-2}}dx \right)^{\frac{n-2}{n}}, \quad u
\in C^{\infty}_0(\xO)  \ .
\end{equation}
Such an inequality was first proven in \cite{Maz} in the special case where $\xO$ is the half space.
In \cite{FMT2} it was proven  under the assumption that $\xO$ is a  weakly mean convex domain, that is, it
satisfies in the distributional sense,
\be\la{c1}
- \Delta d(x) \geq 0, ~~~~~~~~~{\mbox in ~~~\xO}  \ .
\ee
 In   \cite{FL}  inequality (\ref{01}) was established  with a
constant $c$ independent of $\xO$, under the stronger hypothesis that  $\xO$ is convex. We note that mean convexity
is equivalent to convexity in $n=2$ dimensions but it is a much  weaker assumption for $n \geq 3$, cf \cite{AK}.

Our interest in this work is in the fractional (non local) Laplacian in a bounded domain. Various
  fractional  $s$--Laplacians ($0<s<1$) have been recently studied, see   \cite{BBC},
  \cite{CT}, \cite{DF},  \cite{FMT1} and references therein.
In  \cite{FMT1} the  limiting case  of obtaining  Hardy--Sobolev--Mazy'a  inequalities for
the half Laplacian   was left open in the case of a domain $\xO$.
 In fact, the half Laplacian is a border line case, since different behaviors
are observed for $s<\frac12$  and  $s >\frac12$. For instance, the fractional Laplacian considered in \cite{D},
satisfies Hardy inequality for  $\frac12<s<1$ but not for $0<s \leq \frac12$, in the case of smooth  bounded domains.
Similarly a  dichotomy appears, for a  different  fractional Laplacian this time,
 in the context of $\Gamma$--convergence of non
local phase transitions, or in the context of non local surface diffusion, see,  \cite{SV},  \cite{CRS}. In our
case certain aymptotics are different for $s>1/2$ than for $s=1/2$ and as a consequence the analysis in \cite{FMT1}
fails for the limiting case of the half Laplacian.

As we have already mentioned, there are several
fractional Laplacians, but in this work   we will focus on the spectral fractional Lapacian that was recently considered
in \cite{CT}. We will do this
as in \cite{FMT1} via   a suitable extension problem in the spirit of \cite{CS}.  In our case  the appropriate
extension problem is the following:
\bea
- \Delta_{(x,y)} u & = & 0,  \hspace{2cm} in ~~  \xO \times (0,\infty) \ ,  \la{01}  \\
 u      & = &     0,      \hspace{2cm}   on~~   \partial \xO \times (0,\infty)      \ ,\la{02}  \\
u(x,0)      & = &   f(x),   \hspace{1,6cm} in ~~  \xO \ ,  \la{03}
\eea
the energy of which is given by
\[
J[u] = \frac12 \dino |\ana_{(x,y)} u(x,y)|^2 \dd  \ .
\]

At this point we recall that the inner radius of a domain $\xO$ is defined as $R_{in} := \sup _{x \in \xO} d(x)$.
We say that the  domain  $\xO$ has finite inner radius whenever  $R_{in} < \infty$.
Our  first result is the following  Trace Hardy-Sobolev-Maz'ya inequality:

\begin{theorem}\la{th11}
\noindent  Let $n \geq 2$  and $\xO \subsetneqq \Ren$  be a  uniformly Lipschitz  domain with finite inner radius
 which  in  addition satisfies
\be\la{c2}
- \Delta d(x) \geq 0, ~~~~~~ {\mbox in ~~~\xO}  \ .
\ee
 Then there exists a positive constant $c$  such that
for all   $u \in C^{\infty}_{0}(\xO \times \R) $ there holds
\be\la{1.A1}
 \dino  |\ana_{(x,y)} u(x,y)|^2 \dd  \geq   \frac{2}{\pi}
 \int_{\xO} \frac{u^2(x,0)}{d(x)}dx + c
 \left(  \int_{\xO} | u(x,0)|^{\frac{2n}{n-1}} dx \right)^{\frac{n-1}{n}}  \ .
\ee
\end{theorem}

We note that the constant $\frac{2}{\pi}$ is the best constant for the corresponding trace Hardy inequality
see Theorem 1 of \cite{FMT1} for the precise statement.

We  will apply Theorem \ref{th11}  to
 the spectral  fractional Laplacian that is defined as follows.
Let $\xO \subset \Ren$ be a bounded  domain, and   $\xl_i$ and  $\phi_i$ be the Dirichlet eigenvalues and
orthonormal  eigenfunctions of the Laplacian,
i.e. $-\Delta \phi_i = \xl_i \phi_i$ in $\xO$, with $\phi_i = 0$ on $\partial \xO$.
 Then,  for $f(x) = \sum c_i \phi_i(x)$ the $s$--fractional Laplacian is defined by
\be\la{1.frc}
(-\Delta)^{s} f(x) =  \sum_{i=1}^{\infty} c_i \xl_i^s  \phi_i(x),~~~~~0<s<1 \ .
\ee
In the sequel we will be interested in the case $s=\frac12$. More precisely, the
 Hardy-Sobolev-Maz'ya inequality  for the  spectral half  Laplacian reads:

\begin{theorem}\la{th12}
\noindent  Let $n \geq 2$  and $\xO \subsetneqq \Ren$  be a  bounded  Lipschitz  domain
 which  in  addition satisfies
\be\la{c3}
- \Delta d(x) \geq 0, ~~~~~~ {\mbox in ~~~\xO}  \ .
\ee
Then there exists a positive constant $c$  such that
for all   $f \in C^{\infty}_{0}(\xO) $ there holds
\be\la{1.A2}
 ((-\Delta)^{\frac12} f , f )_{\xO} \geq \frac{2}{\pi}
 \int_{\xO} \frac{f^2(x)}{d(x)}dx + c
 \left(  \int_{\xO} | f(x)|^{\frac{2n}{n-1}} dx \right)^{\frac{n-1}{n}}  \ .
\ee
\end{theorem}
Again, the constant  $\frac{2}{\pi}$ is the best constant for the corresponding  Hardy inequality
see Theorem 3 of \cite{FMT1} for the precise statement.

We note that the proof is based on the following crucial estimate that was missing in \cite{FMT1}.

\begin{theorem} \la{th13}
 Let  $\xO \subsetneqq \Ren$  be a domain with finite inner radius $R_{in}$ and is
 such that
\[
- \Delta d(x) \geq 0, ~~~~~~ {\mbox in ~~~\xO}  \ .
\]
 If in addition  $A+1>0$, then
 for all
 $u \in C^{\infty}_{0}(\Ren \times \R) $ there holds
\bea\la{1.A3}
\dino \frac{ y^{A+2} d^{A+2}  }{(d^2+y^2)^{A+2}} |\ana u|^2 \dd
 & + & \frac{(A+1)(4A+9)}{4(2A+5)^2}
 \dino \frac{ y^{A+1} d^{A+1}X }{(d^2+y^2)^{A+\frac32}}(- \Delta d) u^2 dxdy     \nonumber   \\
 &    \geq  &
 \frac{(A+1)^2}{8(2A+5)^2}  ~~
 \dino  \frac{ y^A d^{A} X^2}{(d^2+y^2)^{A+1}} u^2   \dd  \ ,
\eea

where $X=X(\frac{d(x)}{R_{in}})$ and $X(t) = (1- \ln t)^{-1}$, $0<t \leq 1$.
\end{theorem}
For Theorem \ref{th13} it is important that the  domain has a  finite inner radius.

Using  Theorem \ref{th13} and  quite similar arguments to the ones leading to the proof of Theorems \ref{th11} and
\ref{th12}
 one can establish HSM--inequalities for the   Dirichlet  half Laplacian defined in
\cite{FMT1}. In particular Theorems 4, 5 and 12  of \cite{FMT1}   are valid for the limiting case $s=1/2$.

In Section 2 we give the proof of Theorem \ref{th13} after  establishing a more general result, where weak  mean
convexity of the domain is not required. In the final Section 3 we give the proofs of Theorems  \ref{th11} and
\ref{th12}.

\setcounter{equation}{0}
\section{The proof of Theorem \ref{th13}}\la{sec2}
In this section we will prove Theorem \ref{th13}. In fact
we will  prove a more general result that does not require  any sign assumption  on the measure $- \Delta d(x)$.

\begin{theorem} \la{th21}
 Let $\xO \subsetneqq \Ren$  be a domain with finite inner radius $R_{in}$. If  $A+1>0$, then
 for all
 $u \in C^{\infty}_{0}(\Ren \times \R) $ there holds
\bea\la{bas1}
\dino \frac{ y^{A+2} d^{A+2}  }{(d^2+y^2)^{A+2}} |\ana u|^2 \dd \hspace{9cm} \nonumber   \\
+  \frac{A+1}{4(2A+5)^2}
\dino   \left[  \frac{ y^{A+2} d^{A+1}X^2}{(d^2+y^2)^{A+2}} +
4(A+2) \frac{ y^{A+2} d^{A+3}X}{(d^2+y^2)^{A+3}} \right] (- \Delta d) u^2 dxdy  \nonumber   \\
\geq
\frac{(A+1)^2}{8(2A+5)^2}  ~~
 \dino  \frac{ y^A d^{A} X^2}{(d^2+y^2)^{A+1}} u^2   \dd  \ ,
\eea
where $X=X(\frac{d(x)}{R_{in}})$ and $X(t) = (1- \ln t)^{-1}$, $0<t \leq 1$.
\end{theorem}
From this estimate we have:

\noindent
{\em Proof of  Theorem \ref{th13}:} The result follows from Theorem \ref{th21} using the
sign assumption  $- \Delta d(x) \geq 0$ in $\xO$.

\finedim

    The rest of this section  is devoted to the proof of Theorem \ref{th21}. We first present some auxilliary
Lemmas. Our first Lemma covers a  limiting case  of Lemma 8  of \cite{FMT1}.
\begin{lemma} \la{lem22}
 Let $\xO \subsetneqq \Ren$ be a domain with finite inner radius $R_{in}$.
If $A$ and  $B$ are constants such that $A+1>0$ and  $B+1>0$  then for all
 $u \in C^{\infty}_{0}(\Ren \times \R) $ there holds
\bea\la{1ab1}
 \dino \frac{ y^A d^{B} X^2}{(d^2+y^2)^{\frac{A+B+2}{2}}} | u| \dd &    \leq&
 \frac{A+B+2}{A+1}  \dino \frac{ y^{A+2} d^{B}X }{(d^2+y^2)^{\frac{A+B+4}{2}}}(- d\Delta d) | u| dxdy
 \nonumber   \\
& &    +  \frac{A+B+3}{A+1}
\dino \frac{ y^{A+1} d^{B} X }{(d^2+y^2)^{\frac{A+B+2}{2}}} |\ana u| \dd \ ,
\eea
where $X=X(\frac{d(x)}{R_{in}})$ and $X(t) = (1- \ln t)^{-1}$, $0<t \leq 1$.
\end{lemma}

\noindent
{\em Proof:} Integrating  by parts in the $y$-variable we compute
\bea\la{ipp22}
(A+1) \dino \frac{ y^A d^{B}X^2}{(d^2+y^2)^{\frac{A+B+2}{2}}} | u| \dd  \leq   (A+B+2)
 \dino \frac{ y^{A+2} d^{B}X^2}{(d^2+y^2)^{\frac{A+B+4}{2}}} | u| \dd   \nonumber    \\
 + \dino \frac{ y^{A+1} d^{B}X^2}{(d^2+y^2)^{\frac{A+B+2}{2}}} | u_y| \dd.
\eea
In the previous calculation there is no boundary term due to our assumptions.
To continue we will estimate the first term in the right hand side above.  To this end we define
the vector field $\vec{F}$ by
\be\la{fab1}
\vec{F}(x,y) := \left( \frac{ y^{A+2} d^{B+1}X \ana d}{(d^2+y^2)^{\frac{A+B+4}{2}}},
 ~ \frac{y^{A+3} d^{B} X}{(d^2+y^2)^{\frac{A+B+4}{2}}} \right).
\ee
We then have
\be\la{id1}
\dino {\rm div} \vec{F} |u| dxdy = - \dino \vec{F} \cdot \ana |u| dx dy \leq  \dino | \vec{F}| |\ana u| dxdy.
\ee
We note that because of our assumptions $A+1>0$ and $B+1>0$, there are no boundary terms in
(\ref{id1}).
Straightforward calculations (in the  sense of distributions)    show that,
\be\la{div1}
{\rm div} \vec{F} = \frac{ y^{A+2} d^{B}X  (d \Delta d)}{(d^2+y^2)^{\frac{A+B+4}{2}}} +
  \frac{ y^{A+2} d^{B}X^2 }{(d^2+y^2)^{\frac{A+B+4}{2}}},
\ee
and
\be\la{mod2}
| \vec{F} |=  \frac{ y^{A+2} d^{B}X }{(d^2+y^2)^{\frac{A+B+3}{2}}}
 \leq \frac{ y^{A+1} d^{B} X }{(d^2+y^2)^{\frac{A+B+2}{2}}}.
\ee
From  (\ref{id1})--(\ref{mod2}) we get
\bea\la{1ab12}
 \dino \frac{ y^{A+2} d^{B}X^2}{(d^2+y^2)^{\frac{A+B+4}{2}}} | u| \dd   \leq
 \dino \frac{ y^{A+2} d^{B}X  (- d \Delta d)}{(d^2+y^2)^{\frac{A+B+4}{2}}}  | u| \dd
\nonumber  \\
+  \dino \frac{ y^{A+1} d^{B} X }{(d^2+y^2)^{\frac{A+B+2}{2}}}  |\ana u| \dd \ .
\eea
Combining (\ref{ipp22}) and (\ref{1ab12}) the result follows.

\finedim

We next obtain the $L^2$--analogue of  Lemma \ref{lem22}:
\begin{lemma} \la{lem23}
 Let $\xO \subsetneqq \Ren$ be a domain with finite inner radius $R_{in}$.
If $A$ and  $B$ are constants such that $A+1>0$ and  $B+1>0$  then for all
 $u \in C^{\infty}_{0}(\Ren \times \R) $ there holds
\bea\la{21a}
 \dino \frac{ y^A d^{B} X^2}{(d^2+y^2)^{\frac{A+B+2}{2}}} u^2 \dd &    \leq&
 \frac{2(A+B+2)}{A+1}  \dino \frac{ y^{A+2} d^{B}X }{(d^2+y^2)^{\frac{A+B+4}{2}}}(- d\Delta d)u^2 dxdy
 \nonumber   \\
& & \hspace{-1cm}   +  \frac{4(A+B+3)^2}{(A+1)^2}
\dino \frac{ y^{A+2} d^{B}  }{(d^2+y^2)^{\frac{A+B+2}{2}}} |\ana u|^2 \dd \ ,
\eea
where $X=X(\frac{d(x)}{R_{in}})$ and $X(t) = (1- \ln t)^{-1}$, $0<t \leq 1$.
\end{lemma}
\noindent
{\em Proof:} We apply Lemma \ref{lem22} to  $u^2$.
 We then  use Young's inequality in the last term of the right hand side:
\[
2 y^{A+1}X |u| |\nabla u| \leq \xe y^A X^2  u^2 + \frac{1}{\xe} y^{A+2} |\nabla u|^2,
\]
with
\[
\xe = \frac{A+1}{2(A+B+3)} \ .
\]
   We omit the details.

\finedim

The following is  a variation  of Lemma 6 of \cite{FMT1},
  in the sense that  no assumption on the sign of $(- \Delta d)$ is required.

\begin{lemma}\la{lem24} Suppose that  $\xO \subsetneqq \Ren$ has finite inner radius.
If $A$, $B$ are constants such that $A+1>0$, $B+1>0$, then for all
 $u \in C^{\infty}_{0}(\Ren \times \R) $ there holds
\bea\la{1log}
(B+1)  \dino \frac{ y^A d^{B}X^2}{(d^2+y^2)^{\frac{A+B+2}{2}}} | u| \dd  &\leq &
(A+B+3)\dino \frac{ y^A d^{B+1}X}{(d^2+y^2)^{\frac{A+B+2}{2} }} |\ana u| \dd \ \nonumber
  \\
&  &  \hspace{-7cm} +\dino   \left[  \frac{ y^A d^{B+1}X^2}{(d^2+y^2)^{\frac{A+B+2}{2}}} +
(A+B+2) \frac{ y^A d^{B+3}X}{(d^2+y^2)^{\frac{A+B+4}{2}}} \right] (- \Delta d) | u| dxdy  \ ,
\eea
where $X=X(\frac{d(x)}{R_{in}})$ and $X(t) = (1- \ln t)^{-1}$, $0<t \leq 1$.
\end{lemma}

\noindent
{\em Proof:} Integrating by parts in the $x$-variables we compute
\bea\la{4l1}
(B+1) \dino \frac{ y^A d^{B}X^2}{(d^2+y^2)^{\frac{A+B+2}{2}}} | u| \dd  + 2
\dino \frac{ y^A d^{B}X^3}{(d^2+y^2)^{\frac{A+B+2}{2}}} | u| \dd   \nonumber   \\
 \leq
 \dino \frac{ y^A d^{B+1} X^2  (- \Delta d)}{(d^2+y^2)^{\frac{A+B+2}{2}}} | u| dxdy
+ (A+B+2) \dino \frac{ y^A d^{B+2}X^2}{(d^2+y^2)^{\frac{A+B+4}{2}}} | u| \dd   \nonumber   \\
    +  \dino \frac{ y^A d^{B+1}X^2}{(d^2+y^2)^{\frac{A+B+2}{2}}} |\ana u| \dd.  ~~~~~~~~~~~~~~~
\eea

In the previous calculation there are  no boundary terms  due to our assumptions.
To continue we will estimate the middle term in the right hand side above. To this end we use
(\ref{id1}) with the following choice of the vector field
  $\vec{F}$:
\[\la{f1log}
\vec{F}(x,y) := \left( \frac{ y^A d^{B+3} X \ana d}{(d^2+y^2)^{\frac{A+B+4}{2}} },
 ~ \frac{y^{A+1} d^{B+2} X}{(d^2+y^2)^{\frac{A+B+4}{2}}} \right).
\]
Straightforward calculations show that
\[ \la{div2log}
{\rm div} \vec{F} = \frac{ y^A d^{B+3} X  ( \Delta d)}{(d^2+y^2)^{\frac{A+B+4}{2} }} +
 \frac{ y^A d^{B+2}X^2 }{(d^2+y^2)^{\frac{A+B+4}{2}}},
\] and
 \[\la{mod2log}
 | \vec{F} |=  \frac{ y^A d^{B+2} X
}{(d^2+y^2)^{\frac{A+B+3}{2} }} \leq \frac{ y^A d^{B+1}X
}{(d^2+y^2)^{\frac{A+B+2}{2}}}.
 \]
We then have that
 \bean
  \dino  \frac{ y^A d^{B+2}X^2 }{(d^2+y^2)^{\frac{A+B+4}{2} }}  |u| \dd  & &   \nonumber \\
 \leq
 \dino \frac{ y^A d^{B+3}X }{(d^2+y^2)^{\frac{A+B+4}{2} }}(- \Delta d) | u| dxdy &  + &
\dino \frac{ y^A d^{B+1}X }{(d^2+y^2)^{\frac{A+B+2}{2}}} |\ana u| \dd.
\eean

Combining the above with  (\ref{4l1}) and the fact that $X \leq 1$ we conclude the proof.

\finedim

We next obtain the $L^2$--analogue of  Lemma \ref{lem24}:
\begin{lemma}\la{lem25}
  Suppose that  $\xO \subsetneqq \Ren$ has finite inner radius.
If $A$, $B$ are constants such that $A+1>0$, $B+1>0$, then for all
 $u \in C^{\infty}_{0}(\Ren \times \R) $ there holds
\bea\la{1logl2}
 \dino \frac{ y^A d^{B}X^2}{(d^2+y^2)^{\frac{A+B+2}{2}}}  u^2 \dd  &\leq &
\frac{4(A+B+3)^2}{(B+1)^2}\dino \frac{ y^A d^{B+2}}{(d^2+y^2)^{\frac{A+B+2}{2} }} |\ana u|^2 \dd \ \nonumber
  \\
&  &  \hspace{-6,5cm} + \frac{2}{B+1} \dino  \left[  \frac{ y^A d^{B+1}X^2}{(d^2+y^2)^{\frac{A+B+2}{2}}} +
(A+B+2) \frac{ y^A d^{B+3}X}{(d^2+y^2)^{\frac{A+B+4}{2}}} \right] (- \Delta d) u^2 dxdy  \ ,
\eea
where $X=X(\frac{d(x)}{R_{in}})$ and $X(t) = (1- \ln t)^{-1}$, $0<t \leq 1$.

\end{lemma}
\noindent
{\em Proof:} We apply Lemma \ref{lem24} to $u^2$
 We then  use Young's inequality in the first term of the right hand side:
\[
2 d^{B+1}X |u| |\nabla u| \leq \xe d^B X^2  u^2 + \frac{1}{\xe} d^{B+2} |\nabla u|^2,
\]
with
\[
\xe = \frac{B+1}{2(A+B+3)} \ .
\]
   We omit the details.

\finedim

We are now ready to give the proof of Theorem \ref{th21}:

\noindent
{\em Proof of Theorem \ref{th21}:} We first apply Lemma \ref{lem23} with  $B=A+2$  to get:
\bea\la{1a2}
 \dino \frac{ y^A d^{A+2}X^2}{(d^2+y^2)^{A+2}} u^2 \dd &   \leq &
 \frac{4(2A+5)^2}{(A+1)^2}
\dino \frac{ y^{A+2} d^{A+2}}{(d^2+y^2)^{A+2}} |\ana u|^2 \dd  \nonumber   \\
& &    +
 \frac{4(A+2)}{A+1}  \dino \frac{ y^{A+2} d^{A+3}X }{(d^2+y^2)^{A+3}}(- \Delta d) u^2 dxdy \ .
\eea
We next use Lemma \ref{lem25} with $A=B+2$, to obtain:
\bea\la{1a22}
 \dino \frac{ y^{B+2} d^{B}X^2}{(d^2+y^2)^{B+2} }  u^2 \dd  &\leq &
\frac{4(2B+5)^2}{(B+1)^2}\dino \frac{ y^{B+2} d^{B+2}}{(d^2+y^2)^{B+2}   } |\ana u|^2 \dd \ \nonumber
  \\
&  &  \hspace{-5,5cm} + \frac{2}{B+1} \dino  \left[  \frac{ y^{B+2} d^{B+1}X^2}{(d^2+y^2)^{B+2}} +
2(B+2) \frac{ y^{B+2} d^{B+3}X}{(d^2+y^2)^{B+3}} \right] (- \Delta d) u^2 dxdy  \ ,
\eea

Replacing $B$ by $A$ in (\ref{1a22}) and adding it to (\ref{1a2}) we conclude the proof.

\finedim

\setcounter{equation}{0}
\section{The proofs of Theorems \ref{th11} and \ref{th12}}\la{sec3}

We first establish the following Hardy--Sobolev estimate,
that will be used in an  essential way in the proof of Theorem 1.1. For the definition of the
 `` uniformly Lipschitz domain'' see for instance \cite{FMT1}.

\begin{theorem}\la{th31}
  Let $n \geq 2$  and $\xO \subsetneqq \Ren$  be a uniformly Lipschitz domain  with finite inner
 radius that in addition satisfies
\be\la{3.2}
- \Delta d(x) \geq 0, ~~~~~{\mbox  in}~~~ \xO  \ .
\ee
 Then there exists a positive  constant $c$ such that   for all   $u \in C^{\infty}_{0}(\xO \times \R) $    there holds
\be\la{3.4}
 \int_{0}^{+\infty} \int_{\xO}  |\ana_{(x,y)} u(x,y)|^2 \dd  \geq \frac{2}{\pi}
 \int_{\xO} \frac{u^2(x,0)}{d(x)}dx + c
 \left(\int_{0}^{+\infty}  \int_{\xO} | u(x,y)|^{\frac{2(n+1)}{n-1}} dx dy \right)^{\frac{n-1}{n+1}}  \ .
\ee

\end{theorem}

\noindent {\em Proof of Theorem \ref{th31}:} We first  recall  inequality (2.11), from \cite{FMT1}, with
$s=1/2$ (and $a=0$),  that is
\bea\la{3.8}
 \dino  |\ana u|^2 \dd    & \geq  &
\frac{2}{\pi}  \int_{\xO}  \frac{ u^2(x,0)}{d(x)} dx
+ \dino  |\ana u - \frac{\ana \phi}{\phi} u|^2 \dd  \nonumber \\
 &&  -
 \dino \frac{\Delta \phi}{\phi} u^2 \dd  \ ,
\eea
where
 $\phi$ is  given by
\be\la{3.9}
\phi(x,y)= A \left( \frac{y}{d} \right), ~~~~~~~~
 y >0,~~~~~ x \in \xO  \ ,
\ee
and  $A$ solves  (2.3), (2.4) in  \cite{FMT1}, that is
\be \la{3.9a}
A(t) = 1-\frac{2}{\pi} \arctan t,  \hspace{2cm} t \geq 0.
\ee

The result will follow after  establishing  the following inequality:
\be\la{3.10}
\dino|\ana u - \frac{\ana \phi}{\phi} u|^2 \dd
- \dino \frac{\Delta \phi}{\phi} u^2 \dd \geq
 c
 \left(\int_{0}^{+\infty}  \int_{\xO} | u(x,y)|^{\frac{2(n+1)}{n-1}} dx dy \right)^{\frac{n-1}{n+1}}  \ .
\ee

To this end we start with  the Sobolev  inequality,
\[
\int_0^{+\infty} \int_{\xO}  |\ana u| \dd \geq c
\left(\int_0^{+\infty} \int_{\xO} |u(x,y)|^{\frac{n+1}{n}}
\dd \right)^{\frac{n}{n+1}} \ ,~~~~~~~ u \in
C^{\infty}_{0}(\xO \times \R) \ ,
\]
with the choice $u = \phi^{\frac{2n}{n-1}} v$. Hence we obtain
\bea\la{3.12} \int_0^{+\infty} \int_{\xO}
\phi^{\frac{2n}{n-1}}  |\ana v| \dd + \frac{2n}{n-1}
\int_0^{+\infty} \int_{\xO}
\phi^{\frac{n+1}{n-1}} |\ana \phi|
|v| \dd     \nonumber         \\
 \geq
 c \left( \int_0^{+\infty} \int_{\xO} |\phi^{\frac{2n}{n-1}} v|^{\frac{n+1}{n}}
\dd \right)^{\frac{n}{n+1}}   \ .
\eea
Next we will control the second term of the LHS using Lemma 7 of \cite{FMT1}. To this end we recall that
 we have the following asymptotics (cf Lemma 2 of  \cite{FMT1}),
\be\la{3.14}
 \phi^{\frac{2n}{n-1}} \sim
 \frac{ d^{\frac{2n}{n-1}}}{(d^2+y^2)^{\frac{n}{n-1}}} \ , \hspace{1cm}
\phi^{\frac{n+1}{n-1}} |\ana \phi|
\sim
\frac{  d^{\frac{n+1}{n-1}}}{(d^2+y^2)^{\frac{n}{n-1}}} \ .
\ee
We then use Lemma  Lemma 7 of \cite{FMT1}   with the choice $A= 0$,  $B = \frac{n+1}{n-1}$ and
 $\xG = \frac{n}{n-1}$  taking into
account that $A+B+2-2\xG =1  >0$,
to obtain the estimate
\bea\la{3.14a}
 \int_0^{+\infty} \int_{\xO}
 \frac{  d^{\frac{n+1}{n-1}}}{(d^2+y^2)^{\frac{n}{n-1}}}
|v| \dd  \leq  c_1   \int_0^{+\infty} \int_{\xO}
 \frac{ d^{\frac{2n}{n-1}}}{(d^2+y^2)^{\frac{n}{n-1}}}  | \nabla v| \dd    \nonumber     \\
+ c_2  \int_0^{+\infty} \int_{\xO}
 \frac{  d^{\frac{2n}{n-1}}}{(d^2+y^2)^{\frac{n}{n-1}}}  | v| \dd   \ .
\eea
 From this and (\ref{3.12}) we have that
\[
\int_0^{+\infty} \int_{\xO}
\phi^{\frac{2n}{n-1}}  |\ana v| \dd  + \int_0^{+\infty}
\int_{\xO}  \phi^{\frac{2n}{n-1}}  | v| \dd
\geq
 c \left( \int_0^{+\infty} \int_{\xO} |\phi^{\frac{2n}{n-1}} v|^{\frac{n+1}{n}}
\dd \right)^{\frac{n}{n+1}}  \ .
\]
To continue we next set  $v = |w|^{\frac{2n}{n-1}}$ and apply
Schwartz inequality in the LHS. After a simplification  we arrive
at:
 \be\la{3.15} \dino \phi^2 |\nabla w|^2 \dd + \dino
\phi^2  w^2  \dd \geq c
 \left( \dino |\phi w|^{\frac{2(n+1)}{n-1}} \right)^{\frac{n-1}{n+1}} \ .
\ee
To conclude the proof of the Theorem we need the following estimate:
\be\la{3.16}
c  \dino  \phi^2  w^2  \dd  \leq \dino\phi^2 |\nabla w|^2 \dd  -
 \dino (\Delta \phi) \phi w^2 \dd  \ .
\ee
It is worth noticing that the estimates of \cite{FMT1} that work for $\frac12<s<1$ fail to give (\ref{3.16})
 for the
limiting value $s=\frac12$. It is at this point that we will use the more refined estimate of
Theorem \ref{th13} with $A=0$, that is,
\[
  \dino \frac{  X^2w^2 }{d^2+y^2}  \dd \leq
200 \dino \frac{ y^2 d^2  | \nabla w|^2 }{(d^2+y^2)^2} \dd  +
 18 \dino \frac{ yd X  (-\Delta d)   w^2 }{(d^2+y^2)^{\frac32}}  \dd \ ,
\]
here $X= X\left(\frac{d}{R_{in}}\right)$. This implies
\be \la{3.18}
\frac{1}{R_{in}^2} \dino \frac{  d^2 w^2 }{d^2+y^2}  \dd \leq
 200 \dino \frac{  d^2 | \nabla w|^2  }{d^2+y^2} \dd
 + 18 \dino \frac{ y d (-\Delta d) w^2  }{(d^2+y^2)^{\frac{3}{2}}}  \dd \ .
\ee
Taking into account the asymptotics of $\phi$, cf (\ref{3.14}), as well as the fact that
\[
- \Delta \phi = \frac{2}{\pi} \frac{y}{d^2+y^2} (-\Delta d) \ ,
\]
estimate (\ref{3.18}) leads  to (\ref{3.16}). We omit further details.

\finedim

We next give the proof of  Theorem \ref{th11}.

\noindent {\em Proof of Theorem \ref{th11}:} We will use (\ref{3.8}) where
$\phi$ is given, as before,  by (\ref{3.9}), (\ref{3.9a}).
 The result then will follow  once  we  establish:
\be\la{3.46}
\dino|\ana u - \frac{\ana \phi}{\phi} u|^2 \dd
- \dino \frac{\Delta \phi}{\phi} u^2 \dd \geq
 c
 \left( \int_{\xO} | u(x,0)|^{\frac{2n}{n-1}} dx  \right)^{\frac{n-1}{n}}  \ .
\ee
To this end we start  with  the trace  inequality,
\[
\int_0^{+\infty} \int_{\xO} |\ana u| \dd \geq c
 \int_{\xO} |u(x,0)| dx
 \ ,
\]
 valid for   $u \in C^{\infty}_{0}(\xO \times \R)$. We apply this to  $u =
\phi^{\frac{2n}{n-1}} v$. Hence we obtain
 \bea\la{3.48}
\int_0^{+\infty} \int_{\xO}
\phi^{\frac{2n}{n-1}}  |\ana v| \dd + \frac{2n}{n-1}
\int_0^{+\infty} \int_{\xO}
\phi^{\frac{n+1}{n-1}} |\ana \phi|
|v| \dd
 \geq
 c  \int_{\xO} |\phi^{\frac{2n}{n-1}} v|
dx \ .
\eea
Next we will control the second term of the LHS exactly as we did in (\ref{3.14a}) in
 the proof of Theorem \ref{th31}, to arrive at
\[
\int_0^{+\infty} \int_{\xO}
\phi^{\frac{2n}{n-1}}  |\ana v| \dd  + \int_0^{+\infty}
\int_{\xO}
\phi^{\frac{2n}{n-1}}  | v| \dd \geq
 c\int_{\xO} |\phi^{\frac{2n}{n-1}} v(x,0)|
dx  \ .
\]
To continue we next set  $v = |w|^{\frac{2n}{n-1}}$ and apply
Schwartz inequality in the LHS to get after elementary
manipulations that \bea\la{3.50}
 \left( \dino |\phi w|^{\frac{2(n+1)}{n-1}} \dd  \right)
 \left[  \dino  \phi^2 |\nabla w|^2 \dd + \dino  \phi^2  w^2  \dd \right]  \nonumber   \\
  \geq c
 \left(  \int_{\xO}  |\phi w(x,0)|^{\frac{2n}{n-1}} dx \right)^{2}  \ .
\eea
At this point we use (\ref{3.15}) and then   inequality (\ref{3.16}) to arrive at
\[
\dino\phi^2 |\nabla w|^2 \dd  -
 \dino (\Delta \phi) \phi w^2 \dd \geq c
 \left(  \int_{\xO}  |\phi w(x,0)|^{\frac{2n}{n-1}} dx \right)^{\frac{n-1}{n}}  \ ,
\]
which is equivalent to (\ref{3.46})  after the substitution  $u=\phi |w|$.
 We
omit further details.

\finedim

\noindent {\em Proof of Theorem \ref{th12}:} The result follows from Theorem \ref{th11} since the
harmonic extension  $u(x,y)$ in    $\xO \times [0,\infty)$   of $f$, that is, the solution of (\ref{01})--(\ref{03}),
has energy that satisfies
\be\la{3.60}
\dino |\ana u|^2 \dd
   =    ((-\Delta)^{\frac12} f , f )_{\xO} \ ,
\ee
see (8.5) of \cite{FMT1}.

\finedim

\medskip

{ \bf Acknowledgment} The authors  were  partially supported by the FP7-REGPOT
project ACMAC: Archimedes Center for Modeling, Analysis and Computations of the
University of Crete.
AT was also  partially supported by ELKE grant, University of Crete.

\end{document}